\documentclass[12pt]{article}
\usepackage{geometry}                
\usepackage{graphicx}
\usepackage{amssymb}
\usepackage{epstopdf}
\usepackage{amsmath}
\usepackage{amssymb}
\usepackage{amsthm, enumerate}
\newtheorem{theorem}{Theorem}[section]

\newtheorem{example}[theorem]{Example}

\newtheorem{conjecture}[theorem]{Conjecture}

\newtheorem{observation}[theorem]{Observation}

\newtheorem{question}[theorem]{Question}

\newcommand{\eps}{\varepsilon}
\renewcommand{\epsilon}{\varepsilon}

\title{Tree containment and degree conditions}
\author{Maya Stein\footnote{The author is  affiliated to 
 Department of Mathematical Engineering of the University of Chile, and to the Center for  Mathematical Modeling, UMI 2807 CNRS. She acknowledges support by  CONICYT + PIA/Apoyo a centros cient\'ificos y tecnol\'ogicos de excelencia con financiamiento Basal, C\'odigo AFB170001,
and by Fondecyt Regular Grant 1183080.} \bigskip \\  University of Chile \\ }
\date{}
\begin{document}

\maketitle

\ \\ 
\begin{abstract}
We survey results and open problems relating degree conditions with tree containment in graphs, random graphs, digraphs and hypergraphs, and their applications in Ramsey theory.
\end{abstract}

\section{Introduction}

A fundamental question in extremal graph theory is how to guarantee certain subgraphs by imposing a global condition on the host graph. Often, this is a condition on the degree sequence. Classical examples include Tur\'an's theorem on containment of a complete subgraph, or Dirac's theorem on containment of a Hamilton cycle.
One of the most intriguing open questions in the area
 is to determine  degree conditions a graph $G$ has to satisfy in order to ensure  it contains a fixed tree~$T$, or more generally, all trees of a fixed size. 
 
 Let us start with an easy observation. A greedy embedding argument yields that for $k\in\mathbb N$, a minimum degree $\delta (G)$ of at least $k$ is enough to ensure that each  tree~$T$ with~$k$ edges is a subgraph of $G$. Note that $T$ is not necessarily an induced subgraph. Also note  that although a copy of each $k$-edge tree is present, these copies need not be disjoint.
For instance, if $|V(G)|=k+1$, we are considering a complete graph and its spanning trees. 
  
Although the minimum degree condition  $\delta (G)\ge k$ is tight (it cannot be lowered to $k-1$), the condition seems quite strong. It might not be necessary that {\it all} vertices of the host graph have large degree. For path containment, there is a famous result relying on  the {\it average degree}:  
Erd\H os and Gallai~\cite{Erdos1959} showed in 1959 that if $G$ has  average degree  $d(G)>k$ then $G$ contains a $k$-edge path.
 Erd\H{o}s and S\'os  conjectured in 1963 (see~\cite{Erdos64}) that this bound   on  the  average degree should in fact guarantee  {\it all} trees with $k$ edges to appear as subgraphs. This conjecture, its variants and generalisations, will be one of  the guiding themes of this survey. 
 
 We discuss the Erd\H os-S\'os conjecture in  Section~\ref{sec:average}
and then turn to related questions.
Namely, various other conditions have been suggested that might ensure the appearance  of all trees of some fixed size.
One well known conjecture in this direction is the  Loebl--Kom\'os--S\'os conjecture from 1995 (see~\cite{Loebl95}). This conjecture replaces the  assumption on the average degree with an assumption on  the median degree. We will discuss the  Loebl--Kom\'os--S\'os  conjecture and related results in Section~\ref{sec:median}. 

More recent conjectures with the same conclusion employ a condition on a combination of the maximum and the minimum degree. The first conjecture in this direction is due to Havet, Reed, Wood and the author~\cite{2k3:2016}. The idea is that a minimum degree below~$k$ may still be sufficient to find all fixed-size trees, as long as we require one vertex of large degree in the host graph. This vertex both caters for a possible large degree vertex in the tree~$T$, and ensures we have enough space for the embedding of all of $T$. See Section~\ref{sec:maxmin} for details. 

If we only wish to condition on the minimum degree of the host graph, with no assumptions on the maximum degree,  and if our minimum degree condition is strictly below $k$, it is clearly necessary to exclude some trees, for instance stars, from our considerations. More precisely, it will make sense to add a restriction on the maximum degree of the trees we wish to find. There is a well-known result of Koml\'os, S\'ark\"ozy and Szemer\'edi~\cite{KSS95} from 1995, which had been conjectured by Bollob\'as~\cite{bollobas1978} in 1978. It states that in large graphs $G$,  a minimum degree slightly above~$\frac {|V(G)|}2$ is sufficient to guarantee all bounded degree spanning trees. This result will be another recurring theme of this survey. We discuss variations of Koml\'os, S\'ark\"ozy and Szemer\'edi's result in Section~\ref{sec:min}.

It has also been considered to require, apart from a minimum degree condition, additional properties in the host graph, for instance expansion (in terms of large girth, excluded subgraphs, or neighbourhood conditions). With expansion, and for bounded degree trees, the degree bounds on the host graph can be lowered.  This naturally leads to considering random graphs as well. A conjecture of  Kahn~\cite{kahn} regarding the threshold for containment of bounded degree spanning trees was recently solved by Montgomery~\cite{montgomery}. For several of the above mentioned extremal results for tree containment (such as Koml\'os, S\'ark\"ozy and Szemer\'edi's result), there are resilience versions for random graphs. Also, randomly perturbed graphs have been considered as host graphs.
Expansion and random graphs will be discussed in  Section~\ref{sec:random}.

Some of the above conjectures  have direct applications in Ramsey theory, giving upper bounds on Ramsey numbers of trees. For most trees, however, these bounds do not seem to be sharp, and it might be that the correct numbers need to take into account the relative size of the partition classes of the tree. A conjecture of Burr~\cite{Burr74} from 1974 for Ramsey number of trees, although asymptotically confirmed for bounded degree trees in~\cite{HLT02}, has turned out to be far from correct, leaving plenty of open questions in this area. To date, not even the two-colour Ramsey number of double stars is understood. We will give an overview of the state of the art of Ramsey theory for trees in Section~\ref{sec:ramsey}.

Tree containment is also being studied for oriented trees in oriented graphs and digraphs. It is not sufficient to simply consider the degree of the underlying graph, so even the case of the tournament as a host graph is interesting. In Section~\ref{sec:directed}, we will first look at  two famous conjectures on tree containment in tournaments from the 1980's, due to Sumner~\cite{RW83} and Burr~\cite{BurrO}, respectively. Then, we will highlight  a recent conjecture from~\cite{ABHLSTR} which attempts to generalise the Erd\H os--S\'os conjecture and Burr's conjecture at the same time. Finally, we turn to results  generalising the theorem of Koml\'os, S\'ark\"ozy and Szemer\'edi  to digraphs, and some more open questions. 
 
Finally, tree containment problems have  been translated to the hypergraph setting. We will describe  this thriving area in Section~\ref{sec:hyper}. We cover three types of trees: Tight trees, expansions of trees and Berge trees. Each of these notions corresponds to the respective notion for hyperpaths (and these are the most commonly studied hyperpath notions). For tight trees, Kalai's conjecture (see~\cite{FranklFuredi87}) is widely regarded as an analogue of the Erd\H os--S\'os conjecture for hypergraphs. There is also a generalisation of Koml\'os, S\'ark\"ozy and Szemer\'edi's theorem to tight hypertrees. For expansions of trees and Berge trees, we will present some Erd\H os--S\'os  type results, phrased in terms of their Tur\'an numbers.
 
For the reader's convenience, we summarise here how 
the survey is organised: Section~\ref{sec:average}: Average degree; Section~\ref{sec:median}: Median degree; Section~\ref{sec:min}: Minimum degree; Section~\ref{sec:maxmin}: Maximum and minimum degree; Section~\ref{sec:random}: Expanders and random graphs; Section~\ref{sec:ramsey}: Ramsey theory; Section~\ref{sec:directed}: Directed graphs;  Section~\ref{sec:hyper}: Hypergraphs.

\section{Average degree}\label{sec:average}

The most prominent conjecture on tree containment is a classical conjecture of Erd\H{o}s and S\'os from 1963 which focuses on the average degree. It appeared for the first time in~\cite{Erdos64}.

	\begin{conjecture}[Erd\H{o}s--S\'os conjecture, see~\cite{Erdos64}]\label{ESconj}
	 Every graph with average degree $d(G)>k-1$ contains every tree with $k$ edges as a subgraph.
		\end{conjecture} 
		
		A different way to state this conjecture would be in terms of the extremal number or {\it Tur\'an number} of trees. Namely,  we define as usual the Tur\'an number $ex(n, H)$ of a graph $H$ to be  the largest number of edges an $n$-vertex graph may have without containing $H$ as a subgraph. Then, Conjecture~\ref{ESconj} states that $$ex(n,T)\le\frac{k-1}2n$$ for any $k$-edge tree $T$.
	
		The Erd\H{o}s--S\'os  conjecture is tight for every $k\in\mathbb N$: If $k$ divides $n$, consider the $n$-vertex graph consisting of the union of $\tfrac nk$ disjoint copies of cliques on $k$ vertices. This graph has average degree  $k-1$ but it does not contain any tree  with $k$ edges since  its connected components are too small. 
One can also consider any other $(k-1)$-regular graph, for instance the complete bipartite graph $K_{k-1,k-1}$, which does not contain the star with $k$ edges. 

A structurally different example is given by a complete graph on $n$ vertices, in which all edges inside a set of $\lfloor \frac k2\rfloor -1$ vertices have been deleted. This graph does not contain any balanced tree on $k$ edges. The graph is not extremal, however, as its average degree is slightly lower than the average degree of the examples from the previous paragraph.

		Before giving an overview of the known results concerning the conjecture, let us  insert here a quick observation on the minimum degree we may assume the host graph from the Erd\H{o}s--S\'os conjecture to have.
Since every graph of average degree greater than $k-1$ has a subgraph of minimum degree at least $\frac k2$ and average degree greater than~$k-1$ (this subgraph can be found by successively deleting vertices of too low degree), one can assume that the host graph from  the Erd\H os--S\'os conjecture has minimum degree at least $\frac k2$.

Similarly, one can argue that if we replaced the condition $d(G)>k-1$ with the condition $d(G)>2k-1$, then a greedy embedding of any $k$-edge tree into an appropriate subgraph of $G$  will succeed, and therefore, such a version of Conjecture~\ref{ESconj} trivially holds. The bound $d(G)>2k-1$ can be lowered to $d(G)> 2\frac k{k+1}(k-1)$~\cite{kalai_bip}.

In the early 1990's Ajtai, Koml\'os, Simonovits and Szemer\'edi announced a proof of the Erd\H os--S\'os conjecture for large graphs. Nevertheless, many particular cases have been settled since then, or, in some cases, earlier.

The  results mainly group into four types. First, the conjecture has been verified for special types of trees. Most prominently, and as we mentioned before, a classical result of Erd\H os and Gallai~\cite{Erdos1959} from 1959 implies that  the Erd\H os--S\'os conjecture  holds for paths. 
The Erd\H os--S\'os conjecture is also true for stars and double stars. Indeed, for stars this is trivial, while for double stars it suffices to establish the existence of an edge between a vertex of degree $\ge k$ and a vertex of degree $\ge \frac k2$ in the host graph. Since we can assume that the minimum degree of the host graph is at least~$\frac k2$, such an edge clearly  exists.
Moreover,  it is easy to see that  the Erd\H os--S\'os conjecture holds for all trees having a vertex adjacent to at least $\frac k2$ leaves.
McLennan~\cite{mclennan} showed the conjecture holds for all trees of diameter at most~$4$. Fan, Hong and Liu~\cite{spiders} recently proved the conjecture for all {\it spiders}, i.e.~for all trees having at most one vertex of degree exceeding~$2$.

Second, the Erd\H os--S\'os conjecture  has been verified for special types of host graphs. Brandt and Dobson~\cite{bradob} proved in 1996 that the Erd\H os--S\'os conjecture is true for graphs with girth at least $5$. Sacl\'e and Wo\'zniak~\cite{sacwoz} improved on this result  showing in 1997 that the Erd\H os--S\'os conjecture holds for  all graphs that do not contain $C_4$, the cycle on $4$ vertices.  The conjecture also holds if we exclude certain complete  bipartite subgraphs in the host graph or its complement~\cite{BalaDob, Dob02, hax}, and if the host graph is bipartite~\cite{kalai_bip}.

Third, there are results building on the relation between $k$ and $n$.
In particular, Conjecture~\ref{ESconj} has  been established for several cases when $k$ is very close to~$n$, the order of the host graph 
(note that the largest possible value of $k$ is $k=n-1$ and then we need to find a spanning tree in an almost  complete host graph).
More precisely,  the conjecture holds if $k+1\le n\le k+4$ (for all $k$), and even for the case $n\le k+c$, where~$c$ is any given constant and $k$ is sufficiently large depending on~$c$ (see~\cite{goerlich2016} and references therein).  Furthermore, it is shown in~\cite{BPS3} that if we additionally assume that $k\ge 10^{6}$, then  the Erd\H os--S\'os conjecture holds for  all graphs~$G$ with  $|V(G)|\le (1+10^{-11})k$ and trees $T$ with  $\Delta (T)\le \frac{\sqrt k}{1000}$.

Finally, there are some recent results building on the regularity method, thus only applying to the case when $k$ is linear in $n$, and $n$ is large.
In 2019,  Rozho\v{n}~\cite{rohzon} and independently, the authors of~\cite{BPS3} (and~\cite{thesisGuido}) gave an approximate version of the Erd\H os--S\'os conjecture for trees with linear maximum degree and large dense host graphs. 

\begin{theorem}$\!\!${\rm\bf~\cite{BPS3, rohzon}}\label{thm:ESapprox}
 For each $\delta>0$ there are $n_0$, $\gamma$ such that for each $k$ and for each $n$-vertex graph $G$  with  $n\ge k\ge \delta n\ge \delta n_0$  the following holds.\\ If $G$ satisfies $d(G)\ge (1+\delta)k$, then $G$ contains every $k$-edge tree $T$ with~$\Delta(T)\le\gamma k$. 
\end{theorem}

In~\cite{BPS3}, this is used to obtain the following sharp version of Conjecture~\ref{ESconj} for large dense host graphs, which unfortunately relies on the tree having constant maximum degree.

\begin{theorem}$\!\!${\rm\bf~\cite{BPS3}}\label{thm:ESconst}
For each $\delta>0$ and $\Delta$ there is $n_0$ such that for each $k$ and for each $n$-vertex graph $G$  with  $n\ge k\ge \delta n\ge \delta n_0$  the following holds.\\ If $G$ satisfies $d(G)> k-1$, then $G$ contains every $k$-edge tree $T$ with~$\Delta(T)\le\Delta$.
\end{theorem}

\section{Median degree}\label{sec:median}

A well-known variant of the Erd\H os--S\'os conjecture, which replaces the assumption on the average degree with an assumption on the median degree, is
the Loebl-Koml\'os-S\'os conjecture from 1995. Two variants of this conjecture first appeared in~\cite{Loebl95}.

\begin{conjecture}[($\frac n2$--$\frac n2$--$\frac n2$)--Conjecture~\cite{Loebl95}]\label{conj:nhalf}
Every $n$-vertex graph having at least $\frac n2$ vertices of degree at least $\frac n2$  contains each tree on at most $\frac n2$ vertices as a subgraph.
\end{conjecture}

The ($\frac n2$--$\frac n2$--$\frac n2$)-Conjecture has been attributed to Loebl, while according to~\cite{Loebl95}, Koml\'os and S\'os are the originators of the following variation.

\begin{conjecture}[Koml\'os--S\'os Conjecture~\cite{Loebl95}]\label{conj:KS}
Every $n$-vertex graph having more than~$\frac n2$ vertices of degree at least $k$  contains each tree with $k$ edges as a subgraph.
\end{conjecture}

The following amalgamation came to be called the Loebl--Koml\'os--S\'os Conjecture.

\begin{conjecture}[Loebl--Koml\'os--S\'os Conjecture]\label{conj:LKS}
Every $n$-vertex graph having at least~$\frac n2$ vertices of degree at least $k$  contains each tree with $k$ edges as a subgraph.
\end{conjecture}

Note that the Loebl-Koml\'os-S\'os conjecture neither implies nor is implied by Conjecture~\ref{ESconj}.  

Also note that the bound on the degrees in the conjecture  cannot be lowered, because
because we might need to embed a star. Another example is the disjoint union of cliques of order~$k$, which contains no tree with $k$ edges.

As for the number of vertices of large degree, we do not know of any example making the bound $\frac n2$ sharp. The best example we know of is the following.
If $k$ is odd, consider the complete graph on~$k+1$ vertices and delete all edges inside a set of  $\frac{k+3}2$ vertices. It is easy to
check that this graph has $\frac{k+1}2$ vertices of degree $k$, and it  does not contain the $k$-edge path. Taking the disjoint union of 
several such graphs we obtain examples for other values of $n$ (and we can add a disjoint small  graph to reach any value of $n$). The total number of vertices of large degree is somewhat lower than $\frac n2$, and 
in~\cite{LKS1}, it was conjectured that the number given by this example might be the correct number.

\begin{conjecture}$\!\!${\rm\bf~\cite{LKS1}}
\label{conj:LKSstronger}
Let   $G$ be a graph on $n$ vertices having more than $\frac n2-\lfloor \frac{n}{k+1}\rfloor - (n\mod k+1)$ vertices of degree at least $k$. Then $G$ contains each $k$-edge tree.
\end{conjecture}


The Loebl-Koml\'os-S\'os conjecture
(Conjecture~\ref{conj:LKS}) clearly holds for stars, and it has been proved  for several other special classes of trees. One of the first results of this type is due to
Bazgan, Li, and Wo{\'z}niak~\cite{BLW00}, who  proved the
conjecture for paths in 2000.  Piguet and the author~\cite{PS07} proved that 
Conjecture~\ref{conj:LKS} is true for trees of
diameter at most 5, which improved earlier results of Barr and Johansson~\cite{Barr} and Sun~\cite{Sun07} for smaller diameter.

Conjecture~\ref{conj:LKS} has also been proved for special classes of host graphs.
Soffer~\cite{Sof00} showed that the conjecture is true if the
host graph has girth at least 7. Dobson~\cite{Dob02} proved the
conjecture for host graphs whose complement does not contain the complete bipartite graph
$K_{2,3}$. 

The use of a different approach to the Loebl-Koml\'os-S\'os conjecture based on the regularity method has been initiated  by Ajtai, Koml\'os, and Szemer\'edi~\cite{AKS} in 1995 who solved an approximate version of  Conjecture~\ref{conj:nhalf} for large graphs. 
Their strategy (see also~\cite{KSS95} which appeared around the same time) relies on the regularity method, and has been replicated in similar forms in numerous articles on tree embeddings in large dense graphs. The leading idea in~\cite{AKS} is to cut up the tree into many tiny trees connected by a constant number of  vertices (of possibly very large degree), and additionally, to find a useful matching structure in the regularised host graph. The tiny trees are then embedded into the regulars pairs corresponding to the matching, while the connecting vertices are embedded in suitable clusters that see a large amount of matching edges. 

 Zhao~\cite{Zhao2011} used a refinement of the approach from~\cite{AKS} plus stability arguments to prove the exact version of Conjecture~\ref{conj:nhalf} for large
graphs.
Also using regularity, an approximate version of Conjecture~\ref{conj:LKS}  for $k$ linear in~$n$ was proved by Piguet
and the author~\cite{PS12}. 
 Finally, adding stability arguments, 
 Hladk\'y and Piguet~\cite{HlaPig:LKSdenseExact} and independently, Cooley~\cite{Cooley:2009} succeeded in proving Conjecture~\ref{conj:LKS} for large dense graphs.
   \begin{theorem}$\!\!${\bf\cite{Cooley:2009, HlaPig:LKSdenseExact}}
For every $q>0$ there is  $n_0$ such that for any $n>n_0$ and $k>qn$, each $n$-vertex graph $G$ with at least
$\frac n2$ vertices of degree at least $k$ contains each $k$-edge tree.
\end{theorem}

The regularity method described above fails in sparse host graphs. A new approach covering also this type of host graphs was explored by
 Hladk{\'y}, Koml\'os, Piguet, Simonovits, Szemer{\'e}di and the present author in~\cite{LKS1, LKS2, LKS3, LKS4} (for a 10-page overview of the proof see~\cite{LKSshort}). These authors introduced a  decomposition technique for graphs (stemming from previous work of some of the authors on the Erd\H os--S\'os conjecture) whose output resembles the regularity lemma if applied to a dense graph but is also meaningful in the sparse setting. This enabled them to show the following approximate version of Conjecture~\ref{conj:LKS} for large trees.
 
  \begin{theorem}$\!\!${\bf\cite{LKS1, LKS2, LKS3, LKS4}}\label{thm:LKS}
  For every $\eps>0$ there is $k_0$ such that for every $k\ge k_0$, every $n$-vertex graph having at least $(1+\eps)\frac n2$ vertices of degree at least $(1+\eps)k$ contains each $k$-edge tree  as a subgraph.
  \end{theorem}
   

In~\cite{skew}, Klimo\v sov\'a, Piguet, and Rohzo\v n suggest an interesting generalisation of the Loebl-Koml\'os-S\'os conjecture, inspired by a question of  Simonovits. Let us say that a $k$-edge tree is $r$-skew if one of its colour classes has size at most $r(k+1)$.

\begin{conjecture}[Skew LKS conjecture]\label{conj:skewLKS}
Let   $G$ be a graph on $n$ vertices having more than $rn$ vertices of degree at least $k$. Then $G$ contains each $r$-skew $k$-edge tree.
\end{conjecture}

The authors of~\cite{skew} show an approximate version of this conjecture for large dense graphs. 
Conjecture~\ref{conj:skewLKS} has also been verified for paths and trees of diameter at most five~\cite{rohzonMaster}. Examples similar to the ones given earlier in this section show the conjecture would be close to tight.

\section{Minimum degree}\label{sec:min}

As mentioned in the introduction, 
any $n$-vertex graph $G$ with minimum degree at least $k$ contains every tree with $k$ edges. Clearly, the bound on the minimum degree is tight, as we might have to embed a star. Even if we disregard for a moment stars and other trees having vertices of very large degree, it is not possible to lower the bound on the minimum degree of the host graph. In order to see this, it suffices to consider the union of several disjoint copies of $K_k$ which  does not contain any tree with $k$ edges. 
However, the latter example only works if $k$ divides $n$. In particular, it fails if $k>\frac n2$. So one might suspect that for $k>\frac n2$, a lower minimum degree condition could be sufficient to ensure that $G$ contains all $k$-edge trees that have bounded maximum degree.

In this direction, Bollob\'as~\cite{bollobas1978} conjectured in 1978 that any  graph on~$n$ vertices and minimum degree at least $(1+o(1))\frac{n}{2}$ would contain every spanning tree whose maximum degree is bounded by a constant. This conjecture was proved by Koml\'os, S\'ark\"ozy and Szemer\'edi~\cite{KSS95} in 1995, giving one of the earliest applications of the Blow-up lemma. 

\begin{theorem}[Koml\'os, S\'ark\"ozy and Szemer\'edi~\cite{KSS95}]\label{thm:KSSoriginal} For all $\delta>0$ and $\Delta\in\mathbb N$, there is $n_0$ such that such that every graph $G$ on $n\ge n_0$ vertices with 
$\delta(G)\ge (1+\delta)\frac{n}{2}$ contains each $n$-vertex tree $T$ with  
 $\Delta(T)\le \Delta$.
\end{theorem}

Subsequently, each of the two bounds in Theorem~\ref{thm:KSSoriginal} has been improved. 

\begin{theorem}[Csaba, Levitt, Nagy-Gy\"orgy and Szemer\'edi~\cite{CLNS10}]\label{thm:CLNS} For all $\Delta\in\mathbb N$, there are $n_0$ and $c$ such that every graph $G$ on $n\ge n_0$ vertices with 
$\delta(G)\ge \frac{n}{2}+c\log n$
contains each $n$-vertex tree $T$ with  $\Delta(T)\le \Delta$.
\end{theorem}

\begin{theorem}[Koml\'os, S\'ark\"ozy and Szemer\'edi~\cite{KSS2001}]\label{thm:KSS} For all $\delta>0$, there are $n_0$ and $c$ such that such that every graph $G$ on $n\ge n_0$ vertices with 
$\delta(G)\ge (1+\delta)\frac{n}{2}$ contains each $n$-vertex tree $T$ with  
 $\Delta(T)\le c\frac{n}{\log n}$.
\end{theorem}

The bound on the minimum degree in Theorem~\ref{thm:CLNS} is essentially tight (see~\cite{CLNS10}). 
Also the bound  on the maximum degree in Theorem~\ref{thm:KSS} is essentially best possible. This can be seen by considering  the random graph with edge probability $p=0.9$ which a.a.s.~does not contain a forest of stars of order $\frac{n}{\log n}$ (and  thus  also does not contain any tree containing such a forest).

In contrast to the results from earlier sections, the results from~\cite{CLNS10, KSS95, KSS2001} are all for the case when the tree and the host graph have the same order. In view of the examples from the beginning of the section, we know that these results  cannot be generalised to non-spanning trees in host graphs of smaller minimum degree. 
However, if in addition  we require the host graph to have a connected component of size at least $k+1$, then it does at least contain the $k$-edge path $P_k$. This is the core observation behind the Erd\H os--Gallai Theorem. Let us state the observation here for later reference.
\begin{observation}[Erd\H os-Gallai~\cite{Erdos1959}, Dirac (see~\cite{Erdos1959})]\label{pathk/2}
If $\delta (G)\ge \frac k2$, $G$ is connected and $|V(G)|\ge k+1$ then $P_k\subseteq G$.
\end{observation}
 In order to see that this observation is true, note that a variant of of Dirac's theorem~\cite{dirac} states that every 2-connected $n$-vertex graph $G$ has a cycle of length at least $\min\{n, 2\delta (G)\}$. 
 So,  if~$G$ has a $2$-connected component of size at least $k+1$, then this component contains a cycle of length at least $k$, and thus also a $k$-edge path (possibly using one edge that leaves the cycle). Otherwise, we can embed either the middle vertex of the path,  or a vertex adjacent to the middle edge, into any cutvertex~$x$ of  $G$, and then greedily embed the remainder of the path into two components of $G-x$, using the minimum degree of~$G$. 

This argument, however, only seems to work for the case when the tree we are looking for is the path. Already the following tree, which has only one vertex of degree $>2$, cannot be embedded into all large enough connected graphs obeying the minimum degree condition from above. Assume $3$ divides $k$ and consider the tree obtained from identifying the starting vertices of three distinct $\frac k3$-edge paths. This tree is not a subgraph of the graph obtained from adding an edge between two cliques of size $\lceil\frac k2 +1\rceil$.

Still, there is hope: It has been suggested that requiring one large degree vertex in the host graph might remedy the situation. This vertex will at the same time provide the necessary space in the host graph, and cater for a possibly existing  large degree vertex of the tree. See the next section for details.

\section{Maximum and minimum degree}\label{sec:maxmin}

As noted in Section~\ref{sec:average}, we may assume that the host graph from the Erd\H os--S\'os conjecture  has minimum degree at least $\frac k2$, and as we have seen in the previous section, this alone is not enough to force all $k$-edge trees as subgraphs. However,  a graph $H$ of average degree exceeding $k-1$ does not only have a subgraph $H'$ of minimum degree $\ge\frac k2$, but this subgraph $H'$ also maintains the average degree of~$H$ (that is, $d(H')\ge d(H)$). Therefore, $H'$ has a vertex of degree at least $k$. So, in Conjecture~\ref{ESconj}, we may assume the host graph $G$ to obey the following three conditions: $\delta (G)\ge\frac k2$, $d(G)>k-1$,  and $\Delta(G)\ge k$.

Now, the conditions $\delta (G)\ge\frac k2$  and $\Delta(G)\ge k$ alone are not sufficient for guaranteeing all $k$-edge trees as subgraphs. This is because of a variation of the example given in the penultimate paragraph of Section~\ref{sec:min}:
Adding a universal vertex to the disjoint union of two cliques of size $\lceil\frac k2 +1\rceil$, we obtain a graph $G$ satisfying the maximum and minimum degree conditions from above. But the tree obtained from joining any three trees on roughly $\frac k3$ vertices each  to a new vertex (of degree $3$)  is not contained in $G$.

However, if we elevate either the bound on $\delta (G)$ or the bound on $\Delta(G)$ sufficiently, this example ceases to work. So one may suspect there is
a suitable combination of conditions on the  minimum and the maximum degree of the host graph that might replace the 
condition on the average degree  in the Erd\H os--S\'os conjecture. 
In this spirit, Havet, Reed,  Wood and the present author~\cite{2k3:2016} 
 put forward the following conjecture.

\begin{conjecture}[$\frac 23$-conjecture~\cite{2k3:2016}]\label{23conj}
Every graph of minimum degree at least $\lfloor{\frac{2k}{3}}\rfloor$ and maximum degree at least~$k$ contains each $k$-edge tree. 
\end{conjecture}

In~\cite{BPS1}, Besomi, Pavez-Sign\'e and the present author suggested another combination of bounds on the maximum and the minimum degree of the host graph.

\begin{conjecture}[$2k$-$\frac k2$ conjecture~\cite{BPS1}]
\label{2k,k/2}
Every graph of minimum degree at least $\tfrac{k}{2}$ and maximum degree at least $2k$ contains each $k$-edge tree.
\end{conjecture}

Conjectures~\ref{23conj} and~\ref{2k,k/2}  are asymptotically best possible, as we will discuss later in this section.

Each of the conjectures is clearly true for stars and double stars. They also hold for paths, because of Observation~\ref{pathk/2}. 
In~\cite{BPS1}, both Conjectures~\ref{23conj} and~\ref{2k,k/2} were proved in an approximate form for large dense host graphs and  trees whose maximum degree is bounded by $k^{\frac{1}{49}}$ and $k^{\frac{1}{67}}$, respectively. (For trees of constant maximum degree, the maximum degree of the host graph can even be lowered slightly.)

For  Conjecture~\ref{23conj}, more is known. Recently, Reed and the present author~\cite{RS19a, RS19b} showed  that the conjecture holds  if   we are looking for a spanning tree in a large graph. 
\begin{theorem}$\!\!${\rm\bf\cite{RS19a, RS19b}}
There is $n_0$ such that for every $n\ge n_0$, every $n$-vertex graph of minimum degree at least $\lfloor{\frac{2(n-1)}{3}}\rfloor$ and maximum degree at least~$n-1$ contains each $n$-vertex tree. 
\end{theorem}
This theorem can also be seen as an extension of Theorem~\ref{thm:KSS}: By elevating the bound on the minimum degree of the host graph, we can dispose of the bound on the maximum degree of the tree.

Moreover,
 in~\cite{2k3:2016}, Havet, Reed,  Wood and the present author prove the following two variants of their conjecture.

\begin{theorem}$\!\!${\rm\bf\cite{2k3:2016}}
\label{2:3maint1}  There are a function $f:\mathbb N\to \mathbb N$ and a constant $\gamma>0$ such that  if  for a graph $G$, either of the following holds
\begin{enumerate}[(i)]
\item $\Delta(G)> f(k)$ and $\delta(G)\ge \lfloor \frac{2m}{3}\rfloor$; or
\item  $\Delta(G)\ge k$ and $\delta(G)\ge (1-\gamma)m$,
\end{enumerate} 
then $G$ contains each $k$-edge tree.
\end{theorem}

Theorem~\ref{2:3maint1} confirms that, even if the bounds suggested by Conjectures~\ref{23conj} and~\ref{2k,k/2} should be incorrect, the idea behind the conjectures  is not: It is possible to simultaneously bound the maximum and the minimum degree of a graph (with the bound on the minimum degree strictly below $k$) and as a result, guarantee the appearance of each tree of size $k$ as a subgraph. 

The proof of part~(i) of Theorem~\ref{2:3maint1} is relatively easy, and relies on strategically placing into a maximum degree vertex of the host graph a vertex of the tree that cuts the tree into conveniently sized components. 
It would be very interesting to find an extension of Theorem~\ref{2:3maint1}~(i), with the bound on the minimum degree lowered to~$\frac k2$ (ideally), or to some other number strictly smaller than $\lfloor \frac{2m}{3}\rfloor$.

The following conjecture of Besomi, Pavez-Sign\'e and the author~\cite{BPS2}  tries to correlate the bounds on maximum and the minimum degree of the host graph given in Conjectures~\ref{23conj} and~\ref{2k,k/2}.

\begin{conjecture}[Intermediate range conjecture~\cite{BPS2}]
\label{conj:ell}
For each $\alpha\in[0,\frac 13)$
every graph of minimum degree at least $(1+\alpha)\frac k2$ and maximum degree at least $2(1-\alpha)k$ contains each $k$-edge tree.
\end{conjecture}

As for the earlier conjectures from this section, it can be seen that Conjecture~\ref{conj:ell} holds for stars, double stars and paths. Moreover,
 an approximate version for large dense host graphs and  trees of bounded maximum degree is shown  in~\cite{BPS2}.

Conjecture~\ref{conj:ell} is best possible for certain values of $\alpha$. Namely, it is shown in~\cite{BPS2} that
for all odd $\ell\in\mathbb N$ with $\ell\ge 3$, and for all $\gamma>0$ there are $k\in \mathbb N$, a $k$-edge tree~$T$, and a graph $G$ not containing $T$ such that
 $\delta(G)\geq  (1+\frac{1}{\ell} -\gamma)\frac k2$ and $\Delta (G)\geq 2(1-\frac{1}{\ell}-\gamma)k$. 
Let us  give a quick description of the example from~\cite{BPS2}. 

\begin{example}\label{example-ell}
Consider two 
 complete bipartite graphs $H_i=(A_i, B_i)$ with $|A_i|$ slightly below $\frac{\ell -1}{\ell}k$ and $|B_i|$ slightly below $\frac{\ell +1}{\ell}\cdot\frac k2$. Let $G^*$ be obtained by adding a new vertex~$x$ to $H_1\cup H_2$, such that $x$ is connected to all of $A_1\cup A_2$.
 
Then
 the tree $T^*$ formed by $\ell$ stars of order $\frac{k}{\ell}$ and an additional vertex~$v$ connected to the centres of the stars does not embed in $G^*$.   
 \end{example}

Taking $\ell=3$, or letting $\ell$ converge to $\infty$, we obtain  examples that prove the asymptotic sharpness of Conjectures~\ref{23conj} and~\ref{2k,k/2}. For  Conjecture~\ref{23conj}, one or both of the bipartite graphs $H_i$ may be replaced with a conveniently sized complete graph. A structurally different example for the sharpness of Conjecture~\ref{23conj} is obtained by joining a universal vertex to a bipartite graph whose sides are both slightly below~$\frac 23k$, which does not contain the tree $T^*$ from above.

Example~\ref{example-ell} only shows the asymptotic tightness of Conjecture~\ref{conj:ell} for values of~$\alpha$ that are equal to $\frac 1{\ell}$, for some odd $\ell$. We are not aware of similarly good examples for other values of $\alpha$, and perhaps, the degree conditions from Conjecture~\ref{conj:ell} could be lowered for these values, and there might be jumps. One such jump can be observed when the minimum degree bound is close to  $\frac 23k$. We believe that with $\delta(G)\ge \frac 23k$, we would only need to bound $\Delta(G)$ by $k$ (according to Conjecture~\ref{23conj}). However, with the bound $\delta(G)\ge(1-\eps)\frac 23k$, for any $\eps>0$,  it becomes necessary to bound $\Delta(G)$ by almost $\frac 43k$ (because of Example~\ref{example-ell}).

Let us remark that if we exclude host graphs that are very close to the graph $G^*$, a different set of maximum/minimum degree conditions might be sufficient. In~\cite{BPS2}, it is shown that any large enough graph~$G$ with 
$\delta(G)\geq(1+\delta)\frac{k}{2}$ and 
 $\Delta(G)\geq(1+\delta)\frac 43k$ either looks very much like the graph $G^*$ from Example~\ref{example-ell}, or contains all $k$-edge trees $T$ with  $\Delta(T)\le k^{\frac{1}{67}}$.

We close the section with a new conjecture due to Klimo\v sov\'a, Piguet, and Rohzo\v n which appeared in~\cite{rohzon}. Their conjecture combines the essence of Conjecture~\ref{conj:skewLKS} with the spirit of the minimum/maximum degree conjectures from the present section.

\begin{conjecture}\label{KPR}$\!\!${\rm\bf\cite{rohzon}}
Every $n$-vertex graph with 
$\delta(G)\ge \frac{k}{2}$ and at least $\frac{n}{2\sqrt k}$ vertices of degree at least $k$ contains all $k$-edge trees.
\end{conjecture}

Conjecture~\ref{KPR} would be tight~\cite{rohzon}.

\section{Expanders and random graphs}\label{sec:random}

In the results and conjectures we have seen so far, the degrees in the host graph are of the same order as the tree we wish to embed. Examples showed that this is necessary, even if we bound the maximum degree of the tree.
However, adding the assumption that the host graph has some expansion properties changes the situation. Different types of expansion have been considered for this problem, among these are large girth, guarantees of large neighbourhoods of small sets of vertices, and exclusion of dense bipartite subgraphs. 
Also random graphs fall into this category.

An early result for tree containment in expanding graphs
is due to 
Friedman and Pippenger~\cite{FP87} who
 extended P\'osa's rotation-extension technique~\cite{Posa:HamiltonRandom}  from paths to trees  and 
showed that if each set $X\subseteq V(G)$ with $|X|\le 2k-2$ has at least $(\Delta+1)|X|$ neighbours then $G$ contains all $k$-edge trees of maximum degree $\Delta$. 
This result has been generalised in~\cite{BCPS:AlmostSpanningRandom, Haxell:TreeEmbeddings}.

A use of expansion in the form of large girth is the result by Brandt and Dobson~\cite{bradob} we cited in Section~\ref{sec:average}. They showed more generally that every graph of girth at least $5$  satisfying $\delta(G)\ge\frac k2$ contains every $k$-edge tree with $\Delta (T)\le\Delta(G)$.  A generalisation of this was conjectured by
Dobson~\cite{dob}, and, after preliminary results by Haxell and \L uczak~\cite{haxluk}, confirmed by Jiang~\cite{Jiang01}: For any $t\in\mathbb N_+$, every graph $G$ of girth at least $2t+1$   satisfying $\delta(G)\ge\frac kt$ contains every $k$-edge tree with
 $\Delta (T)\le\delta(G)$. This result was greatly improved by Sudakov and Vondr\'ak~\cite{SudVoTree} who 
 used tree-indexed random walks to 
 show that every graph~$G$ of girth at least $2t+1$  and $\delta(G)\ge d$ contains every  tree  with $\Delta (T)>(1-\eps)d$ and  $|V(T)|=c d^k$ (where $c$ is a constant depending on $\eps$). The same authors also show the requirement of large girth may be replaced with forbidding the host graph to contain a complete bipartite graph $K_{s,t}$ for certain $s,t$.

Some of these results directly apply to  random graphs. One natural possibility in  this setting is to replace the degree conditions with probability thresholds. Then, the main problem amounts to  determining the probability threshold $p=p(n)$ for the binomial random graph\footnote{The graph $G(n,p)$ is defined as a probability space on the set of all graphs on $n$ (fixed) vertices where every edge appears with probability $p$, independently, but we also refer to an element of this space as {\it the random graph} $G(n,p)$. The random graph $G(n,p)$ is said to have a property $\mathcal P$ asymptotically almost surely (a.a.s.) if the probability of $G(n,p)$ having $\mathcal P$ tends to $1$ as $n$ tends to infinity.} $G(n,p)$ to contain asymptotically almost surely (a.a.s.) each tree/all trees from a given class $\mathcal T_n$ of trees. Clearly, as the error probabilities for missing individual trees might add up, 
 there is a difference between containing ``each tree'' and ``all trees'' (the latter  is often referred to as \emph{universality}).

Most of the relevant literature for tree containment in random graphs is focused on spanning trees, or almost spanning trees, and the first case to be tackled was the  path. 
 Koml\'os and Szemer\'edi~\cite{KomSze:HamiltonRandom} and  Bollob\'as~\cite{Boll:HamiltonRandom} showed the threshold for spanning paths is $p=(1+o(1))\frac{\log n}{n}$. The lower bound follows immediately from the fact that for smaller values of~$p$, there are a.a.s.~isolated vertices in $G(n,p)$.
Kahn (see~\cite{kahn}) conjectured that the same threshold also applies to bounded degree trees. Namely, he conjectured that, for each $\Delta$, there is  $C$ such that, given any sequence of $n$-vertex trees $T_n$, each with  maximum degree at most~$\Delta$, the random graph $G(n, C \frac{\log n}n)$ a.a.s.~contains a copy of $T$.
After preliminary results due to a number of authors 
(see e.g.~\cite{AKS:LongRandom, AlKrSu:NearlySpanningTrees, FerNenPet, JoKriSa:Expanders, Kri:SpanningRandom}), 
Montgomery~\cite{montgomery} recently solved Kahn's conjecture. He showed more generally the following statement.
\begin{theorem}$\!\!${\rm\bf\cite{montgomery}}
For each $\Delta > 0$, there is a $C > 0$ such that the random graph $G(n, C \frac{\log n}n)$ almost surely contains a copy of every $n$-vertex tree $T$ with  maximum degree at most $\Delta$.
\end{theorem}

Predating~\cite{montgomery}, some results for
 {\it almost} spanning trees appeared.
 Most importantly, Alon, Krivelevich, and Sudakov~\cite{AlKrSu:NearlySpanningTrees} proved that for all $\eps$ and $\Delta$ there is  $C$   such that $G(n,\frac{C}{n})$ a.a.s.\ contains all trees of order $(1-\epsilon)n$ of maximum degree at most~$\Delta$. The value of the  constant~$C$ was improved by Balogh, Csaba, Pei and Samotij~\cite{BCPS:AlmostSpanningRandom} by using the embedding result from~\cite{Haxell:TreeEmbeddings}. Their result, as well as many of the results we cited in the last two paragraphs, also apply to other  types of  expanding graphs (not only random graphs).

Balogh, Csaba and Samotij~\cite{BCS10} showed a result in the spirit of Theorem~\ref{thm:KSSoriginal} for {\it subgraphs} of random graphs. In order to appreciate their result, let us observe that Theorem~\ref{thm:KSSoriginal} can be stated in terms of  {\it local resilience}: If we delete some edges from the complete graph $K_n$, in a way that at each vertex, at least a $(\frac 12 + \delta)$-fraction of its incident edges is preserved, then the resulting graph still contains all spanning trees of maximum degree $\Delta$ (where $\delta$, $n$ and $\Delta$ are as in the theorem).
Now, in~\cite{BCS10} this is translated to random graphs (and {\it almost} spanning trees).

\begin{theorem}$\!\!${\rm\bf\cite{BCS10}} 
For all $\Delta$, $\epsilon$ and $\delta$ there is $C$ such that after deleting any set of edges from $G(n,\frac{C}{n})$ in a way that at each vertex, at least a $(\frac 12+\delta)$-fraction of the original incident edges  are preserved, then the resulting graph a.a.s.\ contains all trees of order $(1-\epsilon)n$ and of maximum degree at most $\Delta$.
\end{theorem}

There are also some  {\it global resilience} results for random graphs (in this type of result, a fraction of the edges is deleted without any restrictions on the number of edges deleted at each vertex).  Balogh, Dudek and Li~\cite{BDL19} proved a  version of the Erd\H os--Gallai theorem for random graphs. Namely, they determine asymptotically the number of edges a subgraph of $G(n,p)$ needs to have in order to guarantee a $k$-edge path, for different ranges of $p$ and $k$. 
Ara\'ujo, Moreira and Pavez-Sign\'e~\cite{AMPS19} show a version of  the Erd\H os--S\'os conjecture for random graphs, and linear sized trees. More precisely, they show 
 that for all $\Delta$, $\epsilon$ and  $t\in (0,1)$ there is  $C$ such that after deletion of at most a $(1-t-\eps)$-fraction of the edges of $G(n,\frac{C}{n})$, the resulting graph still contains w.h.p.~all  trees of order $tn$ and of maximum degree at most $\Delta$. It seems not to be known whether analogues of the results from~\cite{BCS10} and~\cite{AMPS19} for spanning trees exist. 
 
 Finally, there are  some recent results for {\it randomly perturbed} graphs. This model relies on a graph of linear but very small minimum degree, which is `randomly  perturbed' by adding a few random edges to it. More precisely, we consider the union of a graph~$G_\alpha$ of minimum degree at least $\alpha n$ and the random graph $G(n,p)$, on the same set of vertices. The study of bounded degree spanning trees in this model was initiated by Krivelevich, Kwan and Sudakov~\cite{perturb1}. They determined the threshold $p=\frac Cn$ (with $C$ depending on $\alpha$ and $\Delta$) for containment of a single tree of maximum degree $\Delta$, and conjectured the same threshold for the corresponding universality result. This was confirmed by B\"ottcher, Han, Montgomery, Kohayakawa, Parczyk and Person~\cite{perturb2}. Joos and Kim~\cite{jooskim} show a variant of Theorem~\ref{thm:KSSoriginal} for randomly perturbed graphs.

\section{Ramsey numbers}\label{sec:ramsey}

Both Conjecture~\ref{ESconj} (the Erd\H os--S\'os conjecture) and Conjecture~\ref{conj:LKS}  (the Loebl--Koml\'os--S\'os conjecture) have a direct application in Ramsey theory. Let us start with $2$-colour Ramsey numbers. The  {\it (2-colour) Ramsey number} $R(H_1, H_2)$  of a pair of graphs $H_1, H_2$ is the smallest integer $n$ such that every 2-colouring of the edges of~$K_n$ contains a copy of $H_1$ in the first colour, or a copy of $H_2$ in the second colour. Generalising the notion to classes $\mathcal H_1, \mathcal H_2$ of graphs, we write  $R(\mathcal H_1, \mathcal H_2)$ for smallest integer $n$ such that every 2-colouring of the edges of~$K_n$ contains a copy of {\it each} $H_1\in\mathcal H_1$ in the first colour, or a copy of {\it each} $H_2\in\mathcal H_2$ in the second colour. We write short $R(H)$ ($R(\mathcal H)$) for $R(H,H)$ ($R(\mathcal H, \mathcal H)$). 

Some of the earliest results on Ramsey numbers for trees were the following. In 1967, Gerencs\'er and Gy\'arf\'as~\cite{GerencserGyarfas} determined the Ramsey number of two paths.  They showed  that $$R(P_k,P_\ell)=k+\lfloor\frac{\ell+1}2\rfloor$$ for $k$-edge and $\ell$-edge paths
$P_k$ and $P_\ell$ with $k\ge\ell\ge 2$. 
For stars, the Ramsey number is known to be larger. Harary~\cite{Harary:RecentRamsey} observed in 1972 that $$R(K_{1,k},K_{1,\ell})=k+\ell$$ if at least one of $k$, $\ell$ is odd, and $R(K_{1,k},K_{1,\ell})=k+\ell-1$  in the case that $k$ and $\ell$ are both even.

Conjectures~\ref{ESconj}  and~\ref{conj:LKS} can be applied as follows in the Ramsey setting.
Given a $2$-edge-coloured $K_{k+\ell}$, say with colours red and blue, it is easy to see that either the red graph has average degree greater than $k-1$, or the blue graph has average degree greater than $\ell-1$. Also, either the red graph  has median degree at least~$k$, or the blue graph  has median degree at least $\ell$. Therefore, each of the two conjectures  would imply that every $2$-edge-colouring of $K_{k+\ell}$ with colours red and blue contains either \emph{all} $k$-edge trees in red, or \emph{all} $\ell$-edge trees in  blue, and therefore, $R(\mathcal T_k,  \mathcal T_\ell)\le k+\ell$, where $\mathcal T_j$ is the class of all trees with $j$  edges. If $k$ and $\ell$ are both even, the bound we can infer from Conjecture~\ref{ESconj} is even lower:  In that case $R(\mathcal T_k,  \mathcal T_\ell)\le k+\ell-1$. 
Accordingly, and focusing on the case $k=\ell$, Burr and Erd\H os~\cite{BE76} conjectured in 1976  that $R(\mathcal T_k)\le 2k$, and $R(\mathcal T_k)\le 2k-1$ if $k$ is even. The bound $R(\mathcal T_k)\le 2k$ has been confirmed  for large $k$, by Zhao's solution of Conjecture~\ref{conj:nhalf} for large host graphs~\cite{Zhao2011}.

However, the bound $R(T_k, T_\ell)\le k+\ell$ for a $k$-edge tree $T_k$ and an $\ell$-edge tree $T_\ell$ seems to be far from best possible for non-star trees. As noted above,  the Ramsey number for paths  differs significantly from the Ramsey number  for stars with the same number of edges. Note that paths are (almost) completely balanced trees, while stars are the most  unbalanced trees. Believing this difference to be the reason for the variation in their Ramsey numbers,
Burr~\cite{Burr74} put forward the following conjecture in 1974. He suggested that if $T$ is a tree whose bipartition classes have sizes $t_1, t_2$, with $t_2\ge t_1\ge 2$, then the Ramsey number of $T$ is  $$R_B(T):=\max\{2t_1+t_2, 2t_2\}-1.$$ Standard examples show this number would be best possible, and $R_B(T)$ matches the Ramsey numbers for paths from~\cite{GerencserGyarfas}.

 Haxell, \L uczak, and Tingley~\cite{HLT02} confirmed Burr's conjecture asymptotically for   trees with (linearly) bounded maximum degree in 2002. However, already shortly after the conjecture was posed, Grossman, Harary and Klawe~\cite{GHK79} found that it was not true for certain double stars. A {\it double star} $D_{t_1,t_2}$ is a union of two stars $K_{1,t_1-1}$ and $K_{1,t_2-1}$ whose centres are joined by an edge. The examples from~\cite{GHK79} still allowed for the possibility that Burr's conjecture was off only by one, that is, that  the Ramsey number of any tree $T$ would be bounded by  $R_B(T)+1$. The authors of~\cite{GHK79} conjectured this to be the truth for double stars. This has been confirmed for a range\footnote{The current best results are $R(D_{t_1,t_2})\le R_B(D_{t_1,t_2})+1$ if $t_2\ge 3t_1-2$  (obtained using ad hoc arguments~\cite{GHK79}) and $R(D_{t_1,t_2})\le R_B(D_{t_1,t_2})$ if $t_2\le 1.699t_1+1$ (obtained using flag algebras~\cite{NSZdouble}).} of values of $t_1, t_2$. But recently, Norin, Sun and Zhao~\cite{NSZdouble}   disproved the conjecture from~\cite{GHK79} in general by showing that  the numbers $R(D_{t_1,t_2})$ and $R_B(D_{t_1,t_2})$ differ considerably if $t_2$ lies between $\frac 74 t_1+o(t_1)$ and $\frac {105}{41} t_1+o(t_1)$. In particular, for the case $t_2=2t_1$ they find that
$$R(D_{t_1,t_2})\ge4.2t_1-o(t_1)$$ while $$R_B(D_{t_1,t_2})=4t_1-1.$$

The authors of~\cite{NSZdouble} pose the following question.

\begin{question}[Norin, Sun and Zhao~\cite{NSZdouble}]
Is it true that $R(D_{t_1,t_2})=4.2t_1+o(t_1)$ if $t_2=2t_1$?
\end{question}

A question of 
Erd\H os, Faudree, Rousseau and Schelp~\cite{EFRS82}, who, in 1982, asked whether $R(T) = R_B(T)$ for all trees $T$
with colors classes of sizes $t_1$ and $t_2=2t_1$, has also been answered in the negative by the above mentioned results from~\cite{NSZdouble}. The authors of~\cite{NSZdouble} offer the following alternative.

\begin{question}[Norin, Sun and Zhao~\cite{NSZdouble}]
Is it true that $R(T) \le 4.2t_1+o(t_1)$ for all trees $T$
with colors classes of sizes $t_1$ and $2t_1$?
\end{question}

Another natural question in this context seems to be whether there is an exact version of the asymptotic  results of  Haxell, \L uczak, and Tingley~\cite{HLT02} (which would interesting even if we had to restrict the maximum degree of the tree by, say, a constant). A second question is whether their result can be extended to graphs of slightly larger maximum degree. In the main result from~\cite{HLT02}, the bound on the maximum degree of the tree is $\delta |V(T)|$, with $\delta$ depending on the approximation. On the other hand,  the known counter\-examples to Burr's conjecture are all double stars $D$ of maximum degree exceeding $\frac 7{11}|V(D)|$. So, there might be a chance that for some reasonable constant $c\le\frac 7{11}$, Burr's conjecture still holds for all trees $T$ of maximum degree at most $c|V(T)|$.

\begin{question}
Is there a constant $c$ such that $R(T) = R_B(T)$ for all trees $T$
with $\Delta(T)\le c|V(T)|$?
\end{question}

Let us now briefly look at results and questions for multi-colour Ramsey numbers of trees. The $r$-colour Ramsey number $R_r(H)$  of a  graph $H$ is defined as  the smallest integer $n$ such that every $r$-colouring of the edges of $K_n$ contains a monochromatic copy of $H$. 

Multicolour Ramsey numbers for trees have not been studied much. The most studied case is the path $P_k$. The $3$-colour Ramsey number of the
 $k$-edge path $P_k$ has been conjectured to be $2k$ for even $k$ and $2k+1$ for odd $k$ by Faudree and Schelp~\cite{FS}, and this is best possible. This conjecture  has been confirmed for large $k$ by Gy\'arf\'as, Ruszink\'o, S\'ark\"ozy and Szemer\'edi~\cite{GyRSS}.
For more colours, less is known. Constructions based on affine planes show that $R_r(P_k)\ge 
(r-1)k$ if $r-1$ is a prime power. An upper bound on $R_r(P_k)$ can be obtained by applying the Erd\H os--Gallai theorem to the most popular colour in a given $r$-colouring. This yields $R_r(P_k)\le r(k+1)$. Recently, the latter bound has been improved to $(r-\frac 14)(k+1)+ o(k)$ by Davies, Jenssen and Roberts~\cite{DJR17}.

Multicolour Ramsey numbers for $k$-edge stars were determined by Burr and Roberts~\cite{BuRo} in 1973. They showed that $(k - 1)r + 1 \le R_r(K_{1,k}) \le (k -1)r + 2$. The lower bound is tight if and only if both $k$ and $r$ are even.

General bounds for all trees have also been considered.
Erd\H os and Graham~\cite{ErdGra}  observed that an affirmative answer to the following question  would follow from the Erd\H os-S\'os conjecture. (A version for skew trees would also follow from Conjecture~\ref{conj:skewLKS}.)

\begin{question}[Erd\H os and Graham~\cite{ErdGra}]
Is the $r$-colour Ramsey number for a $k$-edge tree $T_k$ equal to $rk+O(1)$?
\end{question}

The authors of~~\cite{ErdGra} observe that $R_r(T_k)$ is bounded from above by $2rk$. This bound can be obtained using a similar argument as for the 2-colour Ramsey number, and the fact that Conjecture~\ref{ESconj} holds if  the average degree bound is replaced by~$2k$, as we observed in Section~\ref{sec:average}.

\section{Directed graphs}\label{sec:directed}

In this section, we will shift our focus from trees and graphs to their oriented versions, that is, oriented trees and digraphs/oriented graphs. An {\it oriented tree (graph)} is a tree (graph) all whose edges have been given a direction. A  {\it digraph} may have (at most) two edges between a pair of vertices, as long as these go in opposite directions. A  {\it tournament} is an oriented complete graph, and a {\it complete digraph} is a digraph having all possible edges.

Let us start with oriented graphs and trees.
Before we turn to possible generalisations of the results in the earlier sections, let us illustrate how the orientations of the edges bring new difficulties.
Just considering the degree in the underlying graph is clearly not enough. Indeed, it is fairly  easy to construct an orientation of a complete graph $K_{k+1}$ that does not contain the $k$-edge star having all its edges directed inwards, thus preventing even the easiest observation for graphs  to carry over to digraphs.

So, let us 
consider the class of tournaments as possible host graphs. One of the first results on oriented trees in tournaments was established by R\'edei~\cite{redei} in~1934. It states that every tournament on $k+1$ vertices contains the {\it directed $k$-edge path}  (i.e.~the  $k$-edge path having all its edges directed in the same direction). More results on oriented paths appeared in, e.g.~\cite{gruenbaum, rosenfeld, tho86}, until in 2000, Havet and Thomass\'e~\cite{HT00b} showed that with three exceptions, all oriented $k$-edge paths appear in any $(k+1)$-vertex tournament. 
Similar results have been shown for some classes of oriented trees with bounded maximum degree (see~\cite{MyNa18} and references therein).

A generalisation of these results for containment of all oriented trees of some fixed size was conjectured by Sumner in the 1980's.
\begin{conjecture}[Sumner, see~\cite{RW83}]\label{sumner}
Every tournament on $2k$ vertices contains every oriented $k$-edge tree.
\end{conjecture}

This conjecture is best possible, which can be seen by considering a {\it $(k-1)$-regular tournament} (that is, a tournament whose vertices each have in- and out-degree $k-1$) on $2k-1$ vertices, which does not contain the $k$-edge star with all edges directed inwards (or outwards). 

Variants of Conjecture~\ref{sumner} replacing  $2k$ with a larger number are known~\cite{elsahili,HT91,  Havet02, HT00a}. The current best bound is $\frac{21}8k$, and was found by Dross and Havet~\cite{DrossHavet}. Havet and Thomass\'e~\cite{HT00a} showed that Conjecture~\ref{sumner} holds for {\it arborescences}, that is, oriented trees having all their edges directed away from (or towards)  a specific vertex. After proving an approximate version~\cite{KMO:ApproximateSumner}, 
K\"uhn, Mycroft and Osthus~\cite{KMO:ProofSumner} confirmed Sumner's conjecture for large $n$, using the regularity method. For oriented trees of bounded degree, the size of the host tournament can be lowered to $n+o(n)$ (see~\cite{KMO:ProofSumner, MyNa18, naiaThesis}).

In 1996, Havet and Thomass\'e proposed that the size of the host tournament can be smaller if we  add a restriction on the number of leaves of the tree. This gives the following generalisation of Conjecture~\ref{sumner}.

\begin{conjecture}[Havet and Thomass\'e, see~\cite{Havet03}]\label{conj:HT}
 Let $T$ be an oriented $k$-edge tree
 with $\ell$ leaves. Then every tournament on $k +\ell$ vertices
contains a copy of $T$.
\end{conjecture}

Note that for oriented stars, Conjecture~\ref{conj:HT} gives the same bound as Conjecture~\ref{sumner}, but for other trees, the bound is lower. As we saw above, if $T$ is a path, the tournament may be by one smaller than required by Conjecture~\ref{conj:HT}.
For  progress on Conjecture~\ref{conj:HT} see~\cite{CeroiHavet, HT91, Havet03}.

Turning now to oriented graphs as possible host graphs for oriented trees, there is  a natural generalisation of the  results and conjectures from above, which involves the chromatic number of a digraph. An oriented graph is {\it $n$-chromatic} if the underlying graph has chromatic number $n$.
The well-known Gallai-Hasse-Roy-Vitaver (GHRV) theorem (see e.g.~\cite{Bang-Jensen2018}) states that every $(k+1)$-chromatic oriented graph contains the  directed path with $k$ edges. As any $n$-vertex tournament is $n$-chromatic, this is a generalisation of R\'edei's theorem mentioned above.

An extension of the GHRV theorem to oriented trees was suggested by Burr~\cite{BurrO} in 1980. His conjecture would imply Sumner's conjecture.

\begin{conjecture}[Burr~\cite{BurrO}]\label{conj:BurrOrient}
 Every $2k$-chromatic oriented graph contains each  oriented $k$-edge tree.
\end{conjecture}

A version of Conjecture~\ref{conj:BurrOrient} for large $k$, replacing oriented trees with oriented paths, and `$2k$-chromatic' with `$k$-chromatic', is attributed in~\cite{Bang-Jensen2018} to Bondy. There is also a generalisation of Conjecture~\ref{conj:HT} in the spirit of Conjecture~\ref{conj:BurrOrient}: Every $(k+\ell)$-chromatic digraph contains each  oriented $k$-edge tree having $\ell$ leaves (this was conjectured in~\cite{HT00a}). 
Furthermore, Naia~\cite{naiaThesis} conjectures that for every oriented $k$-edge tree $T$, the minimum $n$ such that every tournament of order $n$ contains $T$ coincides with the minimum $n$ such that every $n$-chromatic oriented graph contains~$T$. This would imply that Conjecture~\ref{conj:BurrOrient} and Conjecture~\ref{sumner} are equivalent.

Conjecture~\ref{conj:BurrOrient} is only known for some specific classes of oriented paths (see~\cite{ABHLSTR} for references) and  for all oriented stars~\cite{naiaThesis}. 
Burr~\cite{BurrO} showed that Conjecture~\ref{conj:BurrOrient} is true if we replace $2k$ with $k^2$, and Addario-Berry, Havet, Linhares Sales, Thomass\'e and Reed~\cite{ABHLSTR} improved this (roughly by a factor of $2$). A better bound only for oriented graphs with large chromatic number is given in~\cite{naiaThesis}.

The authors of~\cite{ABHLSTR} also propose an interesting conjecture of their own. 
In order to be able to state their conjecture, we need a definition. An {\it antidirected} tree is an oriented tree each of whose vertices either has no incoming edges or no outgoing edges.

\begin{conjecture}[Addario-Berry, Havet, Linhares Sales, Thomass\'e and Reed~\cite{ABHLSTR}]\label{conj:antidir}
 Every  digraph $D$ with more than $(k-1)|V(D)|$ edges contains each  antidirected $k$-edge tree.
\end{conjecture}

This conjecture is best possible because of the $(k-1)$-regular tournament which does not contain the $k$-edge out-star, or alternatively, because of the complete di\-graph on~$k$ vertices which does not contain {\it any} oriented $k$-edge tree. Also note that a version of Conjecture~\ref{conj:antidir} for oriented trees that are not antidirected fails  (even if we made the condition on the number of edges stronger): Consider a large complete bipartite graph $G=(A,B)$, and orient all its edges from $A$ to~$B$. The resulting oriented graph only has antidirected subgraphs.

The authors of~\cite{ABHLSTR}  verify Conjecture~\ref{conj:antidir} for antidirected trees of diameter at most~$3$, and they note it is not difficult to see that Conjecture~\ref{conj:antidir} implies Conjecture~\ref{conj:BurrOrient} (and therefore also Conjecture~\ref{sumner}) for antidirected trees.
(This is because  every $2k$-chromatic graph has a subgraph $H$ of minimum degree at least $2k-1$, and thus  $H$ has more than $(k-1)|V(H)|$ edges.)
They also note that 
if
we restrict  Conjecture~\ref{conj:antidir} to symmetric digraphs (a  digraph is {\it symmetric} if all its edges are bidirected), then the conjecture becomes equivalent to the Erd\H os--S\'os conjecture (Conjecture~\ref{ESconj}). So one can interpret Conjecture~\ref{conj:antidir} as a common generalisation of Conjecture~\ref{ESconj} and Conjecture~\ref{conj:BurrOrient}.

It seems, however, slightly dissatisfying that Conjecture~\ref{conj:antidir} only applies to anti\-directed trees. We have seen above that this is necessary, as there are oriented graphs on $2n$ vertices with  $n^2$ edges that do not even contain  a two-edge directed path.   
In fact, any antidirected host graph with enough edges would serve as an example. Now, in order to avoid these examples, one might try  requiring that the vertices of the host digraph, on average, had {\it both} large enough   in-degree $d^-(v)$  {\it and}  large enough out-degree $d^+(v)$. 

That is, defining the {\it semidegree} of a vertex $v$ as $$d^0(v):=\min\{d^-(v), d^+(v)\},$$  we would require the average of the semidegrees, taken over all vertices $v$ of the digraph $D$, to be larger than $k-1$ (or more generally, to be larger than some function of $k$). Although this would  clearly exclude all antidirected host graphs, it is not  sufficient to guarantee all oriented trees as subdigraphs. In order to see this, just consider an appropriate blow-up\footnote{A {\it blow-up} of a digraph $D$ is obtained by replacing each vertex with an independent set of vertices, and adding all edges from such a set $X$ to a set $Y$, if  $X$ and $Y$ originated from vertices $x$ and $y$ belonging to an edge $\overset{\rightarrow}{xy}$ of $D$.} of  a $(k-1)$-edge directed path. (Observe that this example also shows that a na\"ive extension of the Loebl-Koml\'os-S\'os conjecture that replaces $d(v)$ with 
$d^0(v)$ fails.)

Another possibility is to consider the {\it minimum semidegree} $$\delta^0(D):=\min\{d^0(v):v\in V(D)\}$$ of a digraph $D$. 
Using a greedy embedding argument,  it is clear that  any digraph with $\delta^0(D)\ge k$ must contain each oriented $k$-edge tree. 

This trivial bound can be lowered if $k=n$, and the tree is a path.
Indeed, results from~\cite{DBKMOT15, DBM15} imply that if~$D$ is an $n$-vertex digraph with $\delta^0(D)\ge\frac{n}2$, then $D$ contains every orientation of the path on $n$ vertices, and this is sharp. This might extend to oriented paths of smaller  size.
 A first question in this direction would be whether Observation~\ref{pathk/2} extends to oriented graphs.

\begin{conjecture}\label{q1ori}
Does
every oriented graph $D$ with $\delta^0(D)> \frac k2$ contain each oriented $k$-edge path?
\end{conjecture}

If this conjecture is true, it would be sharp. This can be seen by considering, for even $k$, a blow-up of the directed triangle, replacing each vertex with an independent set of size $\frac k2$. The antidirected path with $k$ edges is not contained in this graph. Moreover, the conjecture is true for directed paths, by a result of Jackson~\cite{jackson}. If we replace the bound on the minimum semidegree with $\frac 34 k$, it holds for antidirected paths~\cite{TM}. If the host graph is a tournament, Conjecture~\ref{q1ori} follows from Conjecture~\ref{sumner}.

An analogous question can be asked for digraphs. Observe that now, we need to require, in addition to the minimum semidegree condition, a lower bound on the size of the largest component (in order to prevent the digraph being the union of complete digraphs of order $\frac k2+2$).

\begin{question}\label{q1}
Does
every digraph $D$ with $\delta^0(D)> \frac k2$ having a component of size at least $k+1$ contain each oriented $k$-edge path?\\ If not, can we lift the bound on the minimum semidegree (to some bound strictly below $k$) so that the question can be answered in the affirmative?
\end{question}

If necessary, one might additionally   require a larger component (or a large  strong component).

Let us shift our attention from oriented paths to oriented bounded degree trees.
Mycroft and Naia~\cite{mycroftNaia} used the minimum semidegree notion to give an extension of Theorem~\ref{thm:KSSoriginal} to digraphs.

\begin{theorem}[Mycroft and Naia~\cite{mycroftNaia}]\label{thm:myna}
For all positive real $\alpha, \Delta$ there exists $n_0$ such that for all $n\ge n_0$ every $n$-vertex digraph $D$ with $\delta^0(D)\ge (\frac 12+\alpha)n$ contains every oriented $n$-vertex  tree of maximum degree at most $\Delta$.
\end{theorem}

In view of their result we feel encouraged to ask whether generalisations to digraphs, using the minimum semidegree notion, of the results and conjectures from Section~\ref{sec:maxmin}  exist.
 In particular, if Conjecture~\ref{q1ori} (Question~\ref{q1}) is true, one might try for  results in the spirit of Theorem~\ref{2:3maint1}  and  Conjectures~\ref{23conj}, ~\ref{2k,k/2} and~\ref{conj:ell}.

\begin{question}
Are there constants $c<1$ and $C$ such that
every oriented graph (digraph) $D$ with $\delta^0(D)\ge ck$ that has a vertex $v$ with $d^0(v)\ge Ck$  contains each oriented $k$-edge tree?
\end{question}

Another possibility is to substitute the semidegree with another degree notion. One natural candidate is the {\it total minimum degree} $\delta_{tot}(D)$, which is defined as the minimum of the sums of the in- and out-degrees of the vertices of the digraph~$D$.
Mycroft and Naia asked the following question~\cite[Problem 4.1]{mycroftNaia}.

\begin{question}$\!\!${\rm\bf\cite{mycroftNaia}}
Does Theorem~\ref{thm:myna} remain true if $\delta^0(D)$ is replaced by $\frac{\delta_{tot}(D)}2$?
\end{question}

If this is true, one could ask for similar variants of the other open questions from this section.

We close the section with a short remark on Ramsey numbers for oriented trees.
There are two natural notions. The {\it oriented Ramsey number} $R^{\rightarrow}_k(T)$ of an oriented tree $T$ is the smallest integer $n$ such that every $k$-coloured tournament on $n$ vertices contains a monochromatic copy of $T$. The {\it directed Ramsey number} $R^{\leftrightarrow}_k(T)$ is defined in the same way, replacing the $k$-coloured tournament with a $k$-coloured  complete digraph. 
Early results using these notion focused on directed paths and two colours~\cite{ChvaDir, GyLe}. An interesting insight gives a recent work of Buci\`c, Letzter and Sudakov~\cite{BuLeSu19} who establish a difference in the order of magnitude of the two numbers, by showing that $R^{\rightarrow}_k(T)=c_k|V(T)|^{k}$ and $R^{\leftrightarrow}_k(T)=c_k|V(T)|^{k-1}$. As observed in~\cite{BuLeSu19}, the former of these two equalities would also follow from Conjecture~\ref{conj:BurrOrient} (Burr's conjecture).

\section{Hypergraphs}\label{sec:hyper}

We will only discuss  $r$-uniform hypergraphs, and call such hypergraphs $r$-graphs for short. 
As one might expect, there is more than one natural generalisation of trees to hypergraphs.
In what follows, we will  discuss  \emph{tight $r$-trees}, \emph{linear $r$-trees}, \emph{$r$-expansions} and \emph{Berge $r$-trees}. We refer to~\cite{GPbook, keevashSurvey} for an overview of more T\'uran type results for hypergraphs.

\subsection{Tight hypertrees}
We start our overview with tight hypergraphs. Call an $r$-graph a \emph{tight $r$-tree} if its edges can be ordered such that except for the first edge,  every edge consists of an $(r-1)$-set  contained in some previous edge, and an entirely new vertex.
Note that for instance, the widely studied tight $r$-paths are examples of tight $r$-trees. (A {\it tight $r$-path} has vertices $v_1, \ldots, v_n$ and edges $\{v_i, \ldots, v_{i+k-1}\}$ for $1\le i\le n-k+1$.)

For $r$-graphs and tight $r$-trees, Kalai proposed in  1984 the following natural generalisation of the Erd\H os-S\'os conjecture (see~\cite{FranklFuredi87}). 

\begin{conjecture}[Kalai's conjecture, see~\cite{FranklFuredi87}]\label{kalaiC}
Let $r\ge 2$ and let $H$ be an $r$-graph on $n$ vertices with more than $\frac{k-1}{r} {n \choose r-1}$ edges.
Then $H$ contains every tight $r$-tree $T$ having $k$ edges.
\end{conjecture}

As already  noted in~\cite{FranklFuredi87}, it follows from constructions using a result of R\"odl~\cite{roedlPackCov} (or alternatively, one can use designs whose existence is guaranteed by Keevash's work~\cite{keevashDesigns}) that Conjecture~\ref{kalaiC} is tight as long as certain divisibility conditions are satisfied.

It is not difficult to observe (see e.g.~\cite[Proposition 5.4]{FurediJiang15}) that  any $n$-vertex  $r$-graph on $n$ vertices with more than $(k-1) {n \choose r-1}$ edges contains every tight $r$-tree with~$k$ edges, which is a factor of $r$ away from the conjectured bound. This bound can be proved as follows: successively delete all edges at $(r-1)$-sets of vertices that lie in few edges until arriving at a subhypergraph of large minimum `codegree'. Then, greedily embed the tree.

Not much is known on Kalai's conjecture in general. Restricting the class of host $r$-graphs,  it is known that the conjecture holds if the host $r$-graph $H$ is  $r$-partite~\cite{kalai_bip}.

Restrictions on the type of tight $r$-trees have led to the following results. In 1987, 
Frankl and F\"uredi~\cite{FranklFuredi87} showed that Conjecture~\ref{kalaiC} holds for all `star-shaped' tight $r$-trees, that is, for all tight $r$-trees  whose first edge  intersects each other edge in $r-1$ vertices. 
F\"uredi, Jiang, Kostochka, Mubayi and Verstra\"ete show in~\cite{FJKMV2019} an asymptotic version of Conjecture~\ref{kalaiC} for a broadened  concept of `star-shaped' (the first $c$ edges have to intersect all other edges  in $r-1$ vertices, for a constant $c$), and  in~\cite{FJKMV18} an exact result for a class of tight $3$-trees.
F\"uredi and Jiang~\cite{FurediJiang15} show
Conjecture~\ref{kalaiC} for special types of tight $r$-trees with many leaves.

On the opposite extreme of the spectrum of tight $r$-trees, there are the tight $r$-paths. 
Improving on results of Patk\'os~\cite{patkos}, F\"uredi, Jiang, Kostochka, Mubayi and Verstra\"ete~\cite{FJKMV17} show that for tight $r$-paths the bound  in Conjecture~\ref{kalaiC} can be replaced by $\frac{k-1}2 {n \choose r-1}$ if $r$ is even, and by a similar bound if $r$ is odd.
Moreover, an asymptotic version of Kalai's conjecture for tight $r$-paths whose order is  linear in  the order $n$ of the  host $r$-graph has been established by Allen, B\"ottcher, Cooley and Mycroft~\cite{ABCM17} for large $n$. The authors of~\cite{ABCM17} remark that they do not believe their result to be best possible, arguing that the constructions and designs from~\cite{keevashDesigns, roedlPackCov} only exist when the order of the host graph is much larger than the order of the tight path.

It seems natural to seek  extensions of  other results for graphs to $r$-graphs and tight $r$-trees. 
This has been done for Theorem~\ref{thm:KSSoriginal}. As usual, for any $r$-graph $H$, let $\delta_{i}(H)$ ($\Delta_i(H)$) denote the minimum (maximum) number of edges any $i$-subset of  $V(H)$ belongs to. With this notation, and using hypergraph regularity,
Pavez-Sign\'e, Quiroz-Camarasa, Sanhueza-Matamala and the author~\cite{PQSS} show a version of Theorem~\ref{thm:KSSoriginal} for hypergraphs. Namely, they show that for any $\gamma, \Delta>0$, every large enough $r$-graph $H$ 
with
 $\delta_{k-1}(H) \ge (\frac 12 + \gamma)|V(H)|$ contains each  $r$-tree  $T$ of the same order obeying $\Delta_1(T) \leq \Delta$.
 
 One might also  ask for generalisations of the  results/conjectures from Section~\ref{sec:maxmin} to tight hypergraph trees.
 
 \begin{question}\label{quest}
Is there a function $f$ such that every $r$-graph $H$ with $\delta_{k-1}(H) \ge \frac k2$ and $\Delta_{k-1}\ge f(k)$ contains each  $k$-edge  tight $r$-tree?
\end{question}

More cautiously, one could replace $\frac k2$ in Question~\ref{quest} with $(1-\gamma)k$, for some fixed~$\gamma>0$. 
Perhaps it is also possible to extend Conjecture~\ref{conj:LKS} to tight $r$-trees.

 \begin{question}
Let $H$ be an $r$-graph such that at least $\frac{{n\choose r-1}}2$ of its $(r-1)$-tuples each belong to at least $k$ edges. Does $H$  contain each  $k$-edge  tight $r$-tree?
\end{question}

\subsection{Expansions of trees and linear paths}

The {\it $r$-expansion} of a tree $T$ is the $r$-uniform hypergraph obtained from $T$ by adding to each edge  $r-2$ new vertices. A {\it linear $r$-tree} is obtained from an edge by subsequently  adding any number of new edges that each contain precisely one of the previous vertices. An $r$-expansion of a path is also called a  {\it linear path}, as it  satisfies the definition of a linear tree.

The {\it Tur\'an number} $ex_r(n, H)$ of an $r$-graph~$H$ is defined (in complete analogy to the Tur\'an number of a graph) as the maximum number of edges a hypergraph  can have if it does not contain $H$.
There is a considerable amount of literature on Tur\'an numbers  of expansions. For an  overview we refer to the  survey of Mubayi and Verstra\"ete~\cite{Mubayi2016}.  One of the important results relevant for this survey is the determination of the Tur\'an number  of the linear  $r$-path  with $k$ edges for fixed $r\ge 4$ and~$k$ and large~$n$ by F\"uredi, Jiang and Seiver~\cite{FJS14} using the delta-system-method. The case $r=3$ was solved by Kostochka, Mubayi and Verstra\"ete~\cite{KMV15shadow1} using an approach based on random sampling.

F\"uredi~\cite{Fu} asymptotically determined the Tur\'an number for $r$-expansions of  trees for $r\ge 4$, and conjectured the corresponding asymptotics for $r=3$; this was confirmed by Kostochka, Mubayi and Verstra\"ete in~\cite{KMV17shadow2}. These results relate the Tur\'an number of an $r$-expansion $T$  with the minimum size $\sigma (T)$ of a crosscut of $T$ (where a {\it crosscut} is a set of vertices met by every edge of $T$ in exactly one vertex). More precisely,  for a fixed $r$-expansion  $T$, and for $r\ge 3$,  the Tur\'an number $ex_r(n,T)$ is asymptotically determined as follows~\cite{Fu, KMV17shadow2}:
 $$ex_r(n,T)=\big(\sigma(T)-1+o(n)\big){n\choose r-1}.$$
That this bound is asymptotically best possible can be seen by considering the $r$-graph consisting of all edges containing exactly one vertex from a fixed set of size $\sigma(T)-1$: This $r$-graph does not contain $T$.  
See also~\cite{FurediJiang15} for some related results.

\subsection{Berge hypertrees}

Other recent activity has focused on Berge $r$-trees. A {\it Berge $r$-tree} is an $r$-graph~$H$ such that there is a tree $T$ (i.e.~an acyclic connected $2$-graph), an injection from $V(T)$ to $V(H)$, and a bijection from $E(T)$ to $E(H)$ such that the images of the endpoints of any edge $e\in E(T)$ are contained in the image of $e$. This definition gives the usual definition of a Berge path if $T$ is a path.

The Tur\'an number for Berge $r$-paths  $BP^{(r)}_k$ was almost completely determined by Gy\H ori, Katona and Lemons~\cite{GKL16}, with the last remaining case solved in~\cite{DGMT16}. The bound is $ex(n,BP^{(r)}_k)\le\frac{n(k-1)}{r+1}$ if $r\ge k$ and $ex(n,BP^{(r)}_k)\le\frac{n}{k}{k\choose r}$ if $r<k$, and extremal $r$-graphs are known.

Results for $k$-edge Berge $r$-trees $BT^{(r)}_k$ have been obtained by Gerbner, Methuku and Palmer~\cite{GMP18} and by 
Gy\H ori, Salia, Tompkins and Zamora~\cite{GSTZ19}. If $r\ge k(k-2)$ and the tree we are looking for is not a star, then the bound for Berge $r$-paths  from the previous paragraph applies~\cite{GSTZ19}. In the case $k>r$ the best known bound is 
$ex(n,BT^{(r)}_k)\le\frac{2(r-1)n}{k}{k\choose r}$, although this can be lowered by a factor of $2(r-1)$, thus reaching  the bound for Berge paths from the previous paragraph, if we assume the Erd\H os--S\'os conjecture holds~\cite{GMP18}. These bounds are sharp under certain divisibility conditions.

It would be interesting to see extensions, both to linear trees and to Berge trees, of the results we have seen in Sections~\ref{sec:min} and~\ref{sec:maxmin} for graphs, that is, results that  use conditions on the minimum (and maximum) degree of the host graph, instead of the average degree (i.e.~number of edges).

\bibliographystyle{acm}
\bibliography{trees}

\end{document}